\documentclass[leqno,a4,11pt]{article}
\usepackage{amsmath,amssymb,amsthm}
\numberwithin{equation}{section}

\newtheorem*{thm}{Main Theorem}
\newtheorem*{ckthm}{Theorem 1 (Cho and Ki [5])}
\newtheorem{lem}{Lemma}
\makeatletter
\def\@seccntformat#1{\csname the#1\endcsname.\quad}
\def\section{\@startsection {section}{1}{\z@}{ 2.3ex plus 2ex minus
 -.2ex}{2.3 ex plus .2ex}{\Large\bf}} 
\def\subsection{\@startsection {subsection}{1}{\z@}{ 2.3ex plus 2ex minus
 -.2ex}{2.3 ex plus .2ex}{\large\bf}} 
\makeatother

\begin{document}
 \font\klsl=cmsl10 at 9pt
\font\fraktur=eufm10

\title{\bf Jacobi operators along the structure flow on real
hypersurfaces in a nonflat complex space form}

\author{U-Hang Ki, Hiroyuki Kurihara and Ryoichi Takagi  }

\date{}
\maketitle

\begin{abstract}
Let $M$ be a real hypersurface of a complex space form with  almost
contact metric structure $(\phi, \xi, \eta, g)$. In this paper, we
study real hypersurfaces in a complex space form whose structure
Jacobi operator $R_\xi=R(\cdot,\xi)\xi$ is $\xi$-parallel. In
particular, we prove that the condition $\nabla_{\xi} R_{\xi}=0$
characterizes the homogeneous real hypersurfaces of type $A$ in a
complex projective space or a complex hyperbolic
space when $R_{\xi}\phi S=S\phi R_{\xi}$ holds on
$M$, where $S$ denotes the Ricci tensor of type (1,1) on $M$.
\end{abstract}


  \setcounter{equation}{0}

\section{Introduction}

Let $(M_n(c), \phi, \tilde{g})$ be a complex $n$-dimensional complex
space form with K\"{a}hler structure $(\phi, \tilde{g})$ of constant
holomorphic sectional curvature $4c$ and let $M$ be an orientable
real hypersurface in $M_n(c)$. Then $M$ has an almost contact metric
structure $(\phi, \xi, \eta, g)$ induced from $(\phi, \tilde{g})$.

It is known that there are no real hypersurface with parallel Ricci
tensors in a nonflat complex space form (see [6], [8]). This result
say that there does not exist locally symmetric real hypersurfaces in
a nonflat complex space form. The structure Jacobi operator
$R_\xi=R(\cdot,\xi)\xi$  has a fundamental role in contact geometry. Cho
and the first author started the study on real hypersurfaces in a complex
space form by using the operater $R_{\xi}$ in [3], [4] and [5]. 
Recently Ortega, P\'{e}rez and Santos [12] have proved that there are no
real hypersurfaces in $P_n \mathbb{C}, n \ge 3$ with parallel
structure Jacobi operator $\nabla R_\xi =0$. More generally, such a
result has been extended by [13] due to them.

Now in this paper, motivated by results mentioned above we consideer
the parallelism of the structure Jacobi operator $R_\xi$ in the
direction of the structure vector field, that is $\nabla_\xi R_\xi
=0$.

In 1970's, the third author [14], [15] classified the homogeneous
real hypersurfaces of $P_n \mathbb{C}$ into six types. On the other
hand, Cecil and Ryan [2] extensively studied a Hopf hypersurface,
which is realized as tubes over certain submanifolds in $P_n
\mathbb{C}$, by using its focal map. By making use of those results
and the mentioned work of Takagi, Kimura [10] proved the local
classification theorem for Hopf hypersurfaces of $P_n \mathbb{C}$
whose all principal curvatures are constant. For the case $H_n
\mathbb{C}$, Berndt [1] proved the classification theorem for Hopf
hypersurfaces whose all principal curvatures are constant. Among the
several types of real hypersurfaces appeared in Takagi's list or
Berndt's list, a particular type of tubes over totally geodesic $P_k
\mathbb{C}$ or $H_k \mathbb{C}~(0 \le k \le n-1)$ adding a
horosphere in $H_n \mathbb{C}$, which is called type $A$, has a
lot of nice geometric properties. For example, Okumura [11](resp.
Montiel and Romero [10]) showed that a real hypersurface in $P_n
\mathbb{C}$ (resp. $H_n \mathbb{C}$) is locally congruent to one of
real hypersurfaces of type $A$ if and only if the Reeb flow $\xi$
is isometric or equivalently the structure operator $\phi$ commutes
with the shape operator $H$.

Among the results related $R_\xi$ we mention the following ones.

\begin{ckthm} Let $M$ be a real hypersurface in a nonflat complex space form
$M_n(c)$ which satisfies $\nabla _\xi R_\xi = 0$ and at the same
time $R_\xi H =HR_\xi$.
Then $M$ is a Hopf hypersurface in $M_n(c)$. Further, $M$ is locally
congruent to one of the following hypersurfaces:

{\setlength{\baselineskip}{0mm}
\begin{enumerate}
\item[$(1)$] In cases that $M_n(c)=P_n\mathbb{C}$ with $\eta(H\xi)\neq
0$,
\begin{enumerate}
\item[$(A_1)$] a geodesic hypersphere of radius $r$, where $0<r<
\pi/2$ and $r \neq \pi/4$;
\item[$(A_2)$] a tube of radius $r$
over a totally geodesic $P_k\mathbb{C}~(1 \leq k \leq n-2)$, where
$0<r< \pi/2$ and $r \neq \pi/4$.
\end{enumerate}
\item[$(2)$] In cases $M_n(c)=H_n\mathbb{C}$,
\begin{enumerate}
\item[$(A_0)$] a horosphere;
\item[$(A_1)$] a geodesic hypersphere or a tube over a complex hyperbolic hyperplane $H_{n-1} \mathbb{C}$;
\item[$(A_2)$] a tube over a totally geodesic $H_k\mathbb{C}~(1 \leq k \leq
n-2)$.
\end{enumerate}
\end{enumerate}  }
\end{ckthm}


In this paper we study a real hypersurface in a nonflat complex
space form $M_n(c)$ which satisfies $\nabla_\xi R_\xi=0$ and at the
same time $R_\xi \phi S=S\phi R_\xi$, where $S$ denotes the Ricci
tensor of the hypersurface. We give another characterization of real
hypersurfaces of type $A$ in $M_n(c)$ by above two conditions. The
main purpose of the present paper is to establish Main Thoerem
stated in section 5. We note that the condition $R_{\xi}\phi S=S\phi
R_{\xi}$ is a much weaker condition. Indeed, every Hopf hypersurface
always satisfies this condition.\par

All manifolds in this paper are assumed to be connected and of class
$C^\infty$ and the real hypersurfaces are supposed to be oriented.\\

  \setcounter{equation}{0}


\section{Preliminaries}
 We denote by $M_n(c),~ c\ne 0$ be a nonflat complex space form
 with the Fubini-Study metric $\tilde{g}$ of constant holomorphic sectional
 curvature $4c$ and Levi-Civita connection $\tilde\nabla$.  For an immersed $(2n-1)$-dimensional Riemannian
 manifold $\tau : M\to M_n(c)$, the Levi-Civita connection $\nabla$ of
 induced metric and the shape operator $H$ of the immersion are
 characterized
\begin{equation*}
\tilde{\nabla}_X Y= \nabla _X Y +g(HX, Y)\nu, \quad \tilde{\nabla}_X \nu
=-HX
\end{equation*}
for any vector fields $X$ and $Y$ on $M$, where $g$ denotes the
Riemannian metric of $M$ induced from $\tilde{g}$. In the sequel the
 indeces $i,j,k,l,\dots$ run over the range $\{1,2,\dots , 2n-1\}$
unless otherwise stated.  For a local orthonormal frame field $\{e_i\}$
of $M$, we denote the dual $1$-forms by $\{\theta_{i}\}$. Then the connection forms
$\theta_{ij}$ are defined by 
\begin{equation*}
     d \theta_i + \sum_j \theta_{ij} \wedge \theta_j =0,
   \quad \theta_{ij} + \theta_{ji} =0.
\end{equation*}
Then we have 
\begin{equation*}
 \nabla_{e_i}{e_j} =\sum_k \theta_{kj}(e_i)e_k =\sum_k \varGamma_{kij} e_k,
\end{equation*}
where we put $\theta_{ij}=\sum_k\varGamma_{ijk}\theta_k$.
The structure tensor $\phi =\sum_i \phi_i e_i$ and the structure vector
$\xi =\sum_i \xi_i e_i$ satisfy 
\begin{align}
  & \displaystyle\sum_k \phi_{ik} \phi_{kj} =\xi_i\xi_j-\delta_{ij},\quad\displaystyle\sum_j \xi_j
  \phi_{ij} =0,\quad  \displaystyle\sum_i \xi_i^2=1,\quad
 \phi_{ij}+\phi_{ji}=0,\notag \\ 
  & d\phi_{ij} =\displaystyle\sum_k(\phi_{ik}\theta_{kj} -\phi_{jk}\theta_{ki} -\xi_i
 h_{jk} \theta_k +\xi_j h_{ik}\theta_k),\\
  & d\xi_i=\displaystyle\sum _j \xi_j\theta_{ji}-\displaystyle\sum_{j,k}\phi_{ji} h_{jk}\theta_k.\notag
\end{align}
We denote the components of the shape operator or the second fundamental
 tensor $H$ of $M$ by $h_{ij}$. The components $h_{ij;k}$ of the
 covariant derivative of $H$ are given by $\sum_k h_{ij;k}\theta_k
 =dh_{ij} -\sum_k h_{ik} \theta_{kj} -\sum_k h_{jk}\theta_{ki}$. Then we
 have the equation of Gauss and Codazzi 
\begin{align}
 &R_{ijkl} =c(\delta_{ik}\delta_{jl} -\delta_{il}\delta_{jk}
+\phi_{ik}\phi_{jl} -\phi_{il}\phi_{jk} +2\phi_{ij}\phi_{kl}) +h_{ik}
h_{jl} -h_{il} h_{jk},\\
 &h_{ij;k} -h_{ik;j} =c(\xi_k \phi_{ij} +\xi_i \phi_{kj} -\xi_j \phi_{ik}
    -\xi_i \phi_{jk}),
\end{align}
respectively.\par From (2.2) the structure Jacobi operator $R_{\xi}
=(\varXi_{ij})$ is given by 
\begin{equation}
\varXi_{ij} =\sum_{k,l} h_{ik}h_{jl} \xi_k \xi_l -\sum_{k,l} h_{ij}
 h_{kl} \xi_k \xi_l +c \xi_i \xi_j -c\delta_{ij},
\end{equation}\par From (2.2) the Ricci tensor $S=(S_{ij})$ is
given by
\begin{equation}
 S_{ij} =(2n+1)c\delta_{ij} -3c\xi_i \xi_j +h h_{ij} -\sum_k h_{ik} h_{kj},
\end{equation} 
where $h=\sum_i h_{ii}$.\par

 First we remark

\begin{lem}
Let $U$ be an open set in $M$ and $F$ a smooth function on $U$. We put
 $dF =\sum_i F_i \theta_i$. Then we have
\begin{equation*}
F_{ij} -F_{ji} =\sum_k F_k\varGamma_{kij} -\sum_k F_k\varGamma_{kji}.
\end{equation*}
\end{lem}
Proof. Taking the exterior derivate of $dF =\sum_i F_i\theta_i$, we have
the formula immediately. \hfill$\Box$\vspace{0.2cm}\par

Now we retake a local orthonormal frame field $e_i$ in such a way that 
(1) $e_1 =\phi\nu$, (2) $e_2$ is in the direction of
$\sum_{i=2}^{2n-1}h_{1i} e_i$ and (3) $e_3 =\phi e_2$. Then we have 
\begin{equation}
\xi_1 =1, \  \xi_i=0 \ (i\geq 2), \  h_{1j} =0 \ (j\geq 3) \mbox{ and
 } \phi_{32}=1. 
\end{equation}
We put $\alpha :=h_{11}, \beta :=h_{12}, \gamma :=h_{22}, \varepsilon
 :=h_{23}$ and $\delta :=h_{33}$.\vspace{0.2cm}\par
\underline{Promise}. Hereafter the indeces
 $p,q,r,s,\dots$ run over the range $\{4,5,\dots , 2n-1\}$ unless
 otherwise stated. \vspace{0.2cm}\par

Since $d\xi _i =0$, we have
\begin{align}
&\theta_{12} =\varepsilon\theta_2 +\delta\theta_3 +\displaystyle\sum_p h_{3p}\theta_p, \notag\\
&\theta_{13} =-\beta\theta_1 -\gamma\theta_2 -\varepsilon\theta_3 -\displaystyle\sum_p h_{2p}\theta_p, \\ 
&\theta_{1p} =\displaystyle\sum_q \phi_{qp}h_{q2}\theta_2 +\displaystyle\sum _q\phi_{qp}
 h_{q3}\theta_3 +\displaystyle\sum_q \phi_{qp}h_{qr}\theta_r. \notag
\end{align}
We put
\begin{equation}
\theta_{23} =\sum_i X_i\theta_i,  \quad \theta_{2p} =\sum_i
 Y_{pi}\theta_i, \quad \theta_{3p} =\sum_i Z_{pi}\theta_i. 
\end{equation}
Then it follows from $d\phi_{2i} =0$ that $Y_{pi} =-\sum_q \phi_{pq} Z_{qi}$ or $Z_{pi}
=\sum_q \phi_{pq} Y_{qi}$. The equations (2.4) and (2.5) are rewritten as
\begin{align}
&\varXi_{ij} =-\alpha h_{ij} +h_{1i} h_{1j} +c\delta_{i1}\delta_{j1}
  -c\delta_{ij},\\
& S_{ij} =h h_{ij} -\sum_k h_{ik} h_{jk} -3c\delta_{i1} \delta_{j1}
  +(2n+1)c\delta_{ij},
\end{align}
respectively.\\

  \setcounter{equation}{0}


\section{Real hypersurfaces satisfying $\nabla_\xi R_\xi=0$}

First we assume that $\nabla_{\xi} R_{\xi} =0$. The components
$\varXi_{ij;k}$ of the covariant derivativation of $R_{\xi}
=(\varXi_{ij})$ is given by 
\begin{equation*}
\displaystyle\sum_k\varXi_{ij;k}\theta_k =d\varXi_{ij}
 -\displaystyle\sum_k\varXi_{kj}\theta_{ki}
 -\displaystyle\sum_k\varXi_{ik}\theta_{kj}.
\end{equation*}

Substituting (2.9) into the above equation we have

\begin{equation}
\begin{split}
\displaystyle\sum_k\varXi_{ij;k}\theta_k&\!\!\!\! = -(d\alpha)h_{ij}
 -\alpha dh_{ij} +(dh_{1i})h_{1j} +h_{1i}(dh_{1j})\\ 
& \ +\alpha\displaystyle\sum_k h_{kj}\theta_{ki} -\alpha
 h_{1j}\theta_{1i} -\beta h_{1j}\theta_{2i} -c\delta_{j1}\theta_{1i}\\
& \ +\alpha\displaystyle\sum_k h_{ik}\theta_{kj} -\alpha h_{1i}\theta_{1j} -\beta
 h_{1i}\theta_{2j} -c\delta_{i1}\theta_{1j}.
\end{split}
\end{equation}

Our assumption $\nabla_{\xi} R_{\xi} =0$ is equivalent to
$\varXi_{ij;1}=0$, which can be stated as follows:
\begin{align}
&\varepsilon =0,\quad\alpha\delta +c =0,\quad h_{3p}=0,\\
&(\beta^2 -\alpha\gamma)_1 -2\alpha\displaystyle\sum_p h_{2p} Y_{p1} =0,\\
&(\beta^2 -\alpha\gamma -c)X_1 +\alpha\displaystyle\sum_p h_{2p} Z_{p1} =0,\\
&(\alpha h_{2p})_1 +\alpha\displaystyle\sum_q h_{pq} Y_{q1} +(\beta^2
 -\alpha\gamma)Y_{p1} -\alpha\displaystyle\sum_q h_{2q}\varGamma_{qp1} =0,\\
&\alpha h_{2p} X_1 -\displaystyle\sum_q (\alpha h_{qp} +c\delta_{pq})Z_{q1} =0,\\
&-(\alpha h_{pq})_1 +\alpha h_{2q} Y_{p1} +\alpha\displaystyle\sum_r h_{rq}\varGamma_{rp1}
 +\alpha h_{2p} Y_{q1} +\alpha\displaystyle\sum_r h_{pr}\varGamma_{rq1} =0.
\end{align}

Hereafter we shall use (3.2) without quoting.\par

Furthermore we assume that $R_{\xi}\phi S =S\phi R_{\xi}$. Under the
assumption $\nabla_{\xi} R_{\xi} =0$, we have the following additional
equations
\begin{align}
&(h\delta -\delta^2 +(2n +1)c)h_{2p} =0,\\
&\tilde{R_{\xi}}\tilde{\phi} A =0,\\
&\tilde{R_{\xi}}\tilde{\phi}\tilde{S}
 =\tilde{S}\tilde{\phi}\tilde{R_{\xi}}.
\end{align}
where $A={}^t(h_{24},h_{25},\ldots ,h_{2,2n-1}),\ \tilde{R_{\xi}}
=(\varXi_{pq}), \ \tilde{\phi} =(\phi_{pq}), \ \tilde{S}=(S_{pq})$.\par

 Now, properly speaking, we should denote the equation (2.3) by, e.g.,
 $(23)_{ijk}$. In this paper we denote it by $(ijk)$ simply. Then we
 have the following equations $(112)$--$(q1p)$.\vspace{0.5cm}\\
$(112)\hspace{33pt}\alpha_2 -\beta_1 =0,$\\[2mm]
$(212)\hspace{33pt}\beta_2-\gamma_1 -2\displaystyle\sum_p h_{2p} Y_{p1}
 =0,$\\[2mm]
$(312)\hspace{33pt}(\alpha-\delta)\gamma -\beta X_2 +(\gamma -\delta)X_1
 -\beta^2 -\displaystyle\sum_p h_{2p}Z_{p1} =-c,$\\[2mm]
$(113)\hspace{33pt}\alpha_3 +3\beta\delta -\alpha\beta +\beta X_1 =0,$\\[2mm]
$(213)\hspace{33pt}\beta_3 -\alpha\delta +\gamma\delta
 +(\gamma-\delta)X_1 -\beta^2 -\displaystyle\sum_p h_{2p}Z_{p1}
 =c,$\\[2mm]
$(313)\hspace{33pt}\beta X_3 +\delta_1 =0,$\\[2mm]
$(223)\hspace{33pt}\gamma_3 -2\beta\delta +2\displaystyle\sum_p h_{2p} Y_{p3}
 +(\gamma -\delta)X_2 -\beta\gamma -\displaystyle\sum_p h_{2p}Z_{p2}
 =0,$\\[2mm]
$(323)\hspace{33pt}\displaystyle\sum_p h_{2p} Z_{p3} -\delta_2 -(\gamma
 -\delta)X_3 =0,$\\[2mm]
$(1p1)\hspace{33pt}\alpha_p +\beta Y_{p1} =0,$\\[2mm]
$(12p)\hspace{33pt}\beta_p +2\displaystyle\sum_{q,r} h_{2q}
 \phi_{rq}h_{rp} +\beta Y_{p2} +\alpha\displaystyle\sum_q\phi_{qp}
 h_{2q} =0,$\\[2mm] 
$(13p)\hspace{33pt}-2\delta h_{2p} +\beta Y_{p3} +\alpha h_{2p} -\beta X_p =0,$
\\[2mm]
$(22p)\hspace{33pt}\gamma_p +2\displaystyle\sum_q h_{2q}Y_{qp} -h_{2p2}
 -\displaystyle\sum_q h_{qp} Y_{q2}
 +\beta\sum_q\phi_{qp} h_{2q} +\gamma Y_{p2} +\sum_q
 h_{2q}\varGamma_{qp2} =0,$\\[2mm]
$(23p)\hspace{33pt}\delta X_p +\beta h_{2p} -\gamma X_p
 +\displaystyle\sum_q h_{2q} Z_{qp} -h_{2p3} -\displaystyle\sum_q h_{qp}
 Y_{q3} +\gamma Y_{p3} +\displaystyle\sum_q h_{2q}\varGamma_{qp3}
 =0,$\\[2mm]
$(33p)\hspace{33pt}\delta_p +h_{2p}X_3 -\displaystyle\sum_q h_{qp}
 Z_{q3} +\delta Z_{p3} =0,$\\[2mm]
$(21p)\hspace{33pt}\beta_p +\displaystyle\sum_{q,r} h_{2q}\phi_{rq}
 h_{rp} -h_{2p1} -\displaystyle\sum_q h_{qp} Y_{q1} +\gamma Y_{p1}
 +\displaystyle\sum_q h_{2q}\varGamma_{qp1} =0,$\\[2mm]
$(31p)\hspace{33pt}-\delta h_{2p} +\alpha h_{2p} -\beta X_p +h_{2p} X_1
 -\displaystyle\sum_q h_{qp} Z_{q1} +\delta Z_{p1} =0,$\\[2mm]
$(32p)\hspace{33pt}\delta X_p +\beta h_{2p} -\gamma X_p
 +\displaystyle\sum_q h_{2q} Z_{qp} +h_{2p} X_2 -\displaystyle\sum_q
 h_{pq} Z_{q2} +\delta Z_{p2} =0,$\\[2mm] 
$(1pq)\hspace{33pt}2\displaystyle\sum_{r,s} h_{rp} \phi_{sr} h_{sq}
 -\alpha\displaystyle\sum_r\phi_{rp} h_{rq}
 +\alpha\displaystyle\sum_r\phi_{rq} h_{rp} -\beta Y_{pq} +\beta Y_{qp}
 =-2c\phi_{pq},$\\[2mm]
$(2pq)\hspace{33pt}h_{2pq} +\displaystyle\sum_r h_{rp}Y_{rq}
 -\beta\displaystyle\sum_r\phi_{rp} h_{rq} -\gamma Y_{pq}
 -\displaystyle\sum_r h_{2r}\varGamma_{rpq} - h_{2qp}$\\[2mm]
\hspace{100pt}$-\displaystyle\sum_r h_{rq}Y_{rp}
 +\beta\displaystyle\sum_r\phi_{rq} h_{rp} +\gamma Y_{qp}
 +\displaystyle\sum_r h_{2r}\varGamma_{rqp} =0,$\\[2mm] 
$(q1p)\hspace{33pt}\displaystyle\sum_r h_{rq}\phi_{sr}h_{sp}
 -\alpha\displaystyle\sum_r\phi_{rq}h_{rq} -\beta Y_{qp} -h_{pq1}$\\[2mm]
\hspace{100pt}$+h_{2q}Y_{p1} +\displaystyle\sum_r h_{rq}\varGamma_{rp1}
 +h_{2p}Y_{q1}+\displaystyle\sum_r h_{rp}\varGamma_{rq1}
 =c\phi_{pq}.$\vspace{0.5cm}\par

Remark. We did not write $(p2q), (3pq), (p3q)$ and $(pqr)$ since we can not use them.\\

  \setcounter{equation}{0}


\section{Formulas and Lemmas}
 \underline{Promise}. In the following, we shall abbreviate the expression ``take
 account of the coefficient of $\theta_i$ in the exterior derivative of
 $\cdots$'' to ``see $\theta_i$ of \ $d$ \  of $\cdots$''.\vspace{0.3cm}\par

In this section we study the crucial case where $\beta\ne 0$.  By (3.6)
and $(31p)$ we have 
\begin{equation}
\beta X_p =(\alpha -\delta)h_{2p}.
\end{equation}
This and $(13p)$ imply that
\begin{equation}
\beta Y_{p3} =\delta h_{2p}.
\end{equation}
The equation (3.9) can be rewrittened as
\begin{equation}
\sum_{q,r} (\alpha h_{pq} +c\delta_{pq})\phi_{qr} h_{r2} =0,
\end{equation}
which, together with (4.2), implies
\begin{equation*}
\beta\sum_{q,r}(h_{pq} -\delta\delta_{pq})Z_{q3} =\delta\sum_{q,r}(h_{pq}
 -\delta\delta_{pq})\phi_{qr}Y_{r3}=0.
\end{equation*}
Hence it follows from $(33p)$ and $(1p1)$ that
\begin{equation}
\delta_p =-h_{2p}X_3 \ \ \mbox{and} \ \ 
\alpha_p =-\beta Y_{p1}.
\end{equation}
Thus by (4.4) and $\alpha_p \delta +\alpha\delta_p =0$ obtained from
(3.2) we have 
\begin{equation}
\beta\delta Y_{p1} =-\alpha h_{2p} X_3,
\end{equation}
and so $\sum_p h_{2p} Z_{p1} =0$. By (4.2), we have
\begin{equation}
\sum_p h_{2p} Z_{p3} =\sum_{p,q} h_{2p} \phi_{pq} Y_{q3}
 =\frac{\delta}{\beta}\sum_{p,q} h_{2p}\phi_{pq} h_{2q}=0.
\end{equation}From (3.6), (4.3) and (4.5) we have
\begin{equation}
h_{2p}X_1 =0.
\end{equation}

Now we shall prove the following key lemma.
\begin{lem}
$H(e_2)\in{\rm span}\{e_1,e_2\}$
\end{lem}
Proof. Suppose that $h_{2p}\ne 0$. Then from (4.7) we have $X_1 =0$. We
 can select the vector $e_4$ so that $h_{24}\ne 0$ and $h_{25}=\cdots
 =h_{2,2n-1} =0$. We put $e_5 :=\phi e_4$ and $\rho :=h_{24}(\ne
 0)$. Note that $\phi_{54} =1$. Then by (4.3) we have 
\begin{equation*}
h_{55} =\delta, \quad h_{p5} =0 \ (p\ne 5).
\end{equation*}
Put $p=5$ in $(32p)$. Then by above equation and (4.1) we have $X_5 =0$
 and so $Z_{45} =0$. Thus we have $Y_{55} =0$. Furthermore, put
 $p=q=5$ in $(q1p)$. Then, since $\varGamma_{551} =Y_{55} =0$, we have
\begin{equation}
\alpha_1 =\delta_1 =0.
\end{equation}
Thus, from $(313)$, $(323)$, (4.6) and $(112)$ we have 
\begin{align}
&X_3 =0,\\
&\alpha_2 =\delta_2 =0,\\
&\beta_1 =0.
\end{align}
By (4.4) and (4.9) we have $\alpha_p =\delta_p =0$. Thus it follows from
 $(1p1)$ that 
\begin{equation}
\alpha_p=\delta_p =Y_{p1} =Z_{p1} =0.
\end{equation}

 Now we put $F=\alpha, \ i=1$ and $j=p$ in Lemma 1. Then, from
 (2.7), (4.8), (4.10) and (4.12) we have
\begin{equation*}
0=\alpha_{1p} -\alpha_{p1} =\sum_k \alpha_k\varGamma_{k1p}
 -\sum_k \alpha_k\varGamma_{kp1}
 =\alpha_3(\varGamma_{31p}-\varGamma_{3p1})=\alpha_3 h_{2p}.
\end{equation*}
Thus we have $\alpha_3=0$. Hence it follows from (4.8), (4.10) and
 (4.12) that $\alpha$ and $\delta$ are constant, which, together with
 $(113)$, imply 
\begin{equation}
\alpha =3\delta.
\end{equation}

On the other hand, seeing $\theta_1\wedge\theta_3$ of \ $d$ \ of $\theta_{23}$, we have
\begin{equation}
X_2 =-2\beta.
\end{equation}
Thus, from $(312)$ and (4.13) we have
\begin{equation}
2\delta\gamma +\beta^2 =-c.
\end{equation}

Seeing $\theta_1$ and $\theta_2$ of \ $d$ \ of (4.15) and taking account of
(4.8), (4.11) and $(212)$, we have
\begin{equation}
\gamma_1 =0, \ \beta_2 =0 \ \mbox{ and } \ \gamma_2 =0.
\end{equation}

Moreover, seeing $\theta_5$ of \ $d$ \ of (4.15), we have
\begin{equation}
\delta\gamma_5 +\beta\beta_5 =0.
\end{equation}\par From (3.5) and (4.12) we have
\begin{equation*}
h_{2p1}-\sum_q h_{2q}\varGamma_{qp1} =0.
\end{equation*}
This, together with $(21p)$ and $(12p)$, implies
\begin{align*}
&\beta_p +\rho h_{5p} =0,\\
&\beta_p +2\rho h_{5p} +\alpha\rho\phi_{4p} +\beta Y_{p2} =0.
\end{align*}
Put $p=4,5,6,\ldots,2n-1$ in above two equations to get
\begin{equation}
\begin{array}{l}
\beta_p=
\left\{\begin{array}{ll}
0 & (p\ne 5)\\
-\rho\delta & (p=5)
\end{array}
\right.,\quad Y_{p2}=
\left\{\begin{array}{ll}
0 & (p\ne 5)\\
\rho(\alpha -\delta)/\beta & (p=5)
\end{array}
\right., \vspace{0.3cm}\\
Z_{p2}=
\left\{\begin{array}{ll}
0 & (p\ne 4)\\
-\rho(\alpha -\delta)/\beta & (p=4)
\end{array}
\right..
\end{array}
\end{equation}
Hence from (4.1), (4.2), (4.17) and (4.18) we have
\begin{equation}
\begin{array}{l}
X_p=
\left\{\begin{array}{ll}
0 & (p\ne 4)\\
\rho(\alpha -\delta)/\beta & (p=4)
\end{array}
\right., \quad
Y_{p3}=
\left\{\begin{array}{ll}
0 & (p\ne 4)\\
-\rho\delta/\beta & (p=4)
\end{array}
\right.,\vspace{0.3cm}\\
Z_{p3}=
\left\{\begin{array}{ll}
0 & (p\ne 5)\\
-\rho\delta/\beta & (p=5)
\end{array}\right., \quad
\gamma_p=
\left\{\begin{array}{ll}
0 & (p\ne 5)\\
-\rho\beta & (p=5)
\end{array}
\right..
\end{array}
\end{equation}\par

Now, by $(213)$, $(223)$, (4.15) and (4.19) we have
\begin{equation}
\begin{split}
&\beta_3 =\beta^2 -\gamma\delta =-\alpha\gamma -c = 3\delta (\delta
 -\gamma),\\
&\gamma_3 =3\beta\gamma -4\rho^2\delta/\beta.
\end{split}
\end{equation}

On the other hand, if we put $F=\beta$ and $\gamma$ in Lemma 1, then 
 from (4.11), (4.12), (4.15), (4.16), (4.18) and (4.19) we have
\begin{equation}
\begin{split}
&\gamma\beta_3 +\rho\beta_5 =0,\\
&\gamma\gamma_3 +\rho\gamma_5 =0.
\end{split}
\end{equation}
Eliminating $\beta_3, \beta_5, \gamma_3, \gamma_5, \rho$ and $\beta$
from (4.17), (4.18), (4.20) and (4.21), we have
\begin{equation*}
4\gamma^2 -6\gamma\delta -c =0.
\end{equation*}
Consequently, $\gamma$ is constant, which contradicts $\gamma_5 =\rho\beta$.
\hfill$\Box$\vspace{0.2cm}\par

Owing to Lemma 2 the matrix $(h_{pq})$ is diagonalizable, that is, for a
suitable choice of a orthonormal frame field $\{e_p\}$ we can set
\begin{equation*}
h_{pq} =\lambda_p \delta_{pq}.
\end{equation*}
Then it is easy to see
\begin{equation}
\begin{split}
&\tilde{R}_{\xi} =-((\alpha\lambda_p +c)\delta_{pq}),\\
&\tilde{S}=(\{h\lambda_p -(\lambda_p)^2 +K\}\delta_{pq}),
\end{split}
\end{equation}
where we put $K=(2n+1)c$.\vspace{0.2cm}\par

 Here we shall sum up all equations obtained from Lemma 2.\par

From (4.1), (4.2) and (4.4) we have
\begin{equation}
X_p = Y_{p1} =Z_{p1} =Y_{p3} =Z_{p3} =0,\quad \alpha_p =\delta_p =0.
\end{equation}
This, together with (3.3) and (3.4), imply
\begin{align}
&(\beta^2 -\alpha\gamma)_1 =0,\\
& (\beta^2 -\alpha\gamma -c)X_1 =0.
\end{align}
Put $p=q$ in (3.7). Then we have
\begin{equation}
(\alpha\lambda_p)_1 =0.
\end{equation}
Moreover, from $(112)$--$(32p)$ we have
\begin{align}
&\alpha_2 -\beta_1 =0,\\
&\beta_2 -\gamma_1 =0,\\
&(\alpha -\delta)\gamma -\beta X_2 +(\gamma -\delta)X_1 -\beta^2 =-c,\\
&\alpha_3 +3\beta\delta -\alpha\beta +\beta X_1 =0,\\
&\beta_3 -\alpha\delta +\gamma\delta +(\gamma -\delta) X_1 -\beta^2 =c,\\
&\delta_1 +\beta X_3 =0,\\
&\gamma_3 -2\beta\delta +(\gamma -\delta)X_2 -\beta\gamma =0,\\
&\delta_2 +(\gamma -\delta)X_3 =0,\\
&\beta_p =0,\\
&Y_{p2} =0,\quad Z_{p2} =0,\\
&\gamma_p =0.
\end{align}
It follows from $(q1p)$ and (3.7) that
\begin{equation}
\alpha\beta Y_{pq} =\alpha \lambda_p\lambda_q \phi_{pq} -\alpha^2
 \lambda_p\phi_{pq} +\alpha_1 \lambda_p\delta_{pq} -c\alpha \phi_{pq}.
\end{equation}
From this and $(2pq)$ we have
\begin{equation}
\beta^2(\lambda_p +\lambda_q)\phi_{pq}
 -(\lambda_p -\gamma)(\lambda_p\lambda_q -\alpha\lambda_q -c)\phi_{pq}
 -(\lambda_q -\gamma)(\lambda_p\lambda_q -\alpha\lambda_p -c)\phi_{pq} =0.
\end{equation}\\

  \setcounter{equation}{0}


\section{Proof of Main theorem}

In this section we prove
\begin{thm}
Let $M$ be a real hypersurface of a complex space form
$M_n(c),~c\neq 0$ which satisfies $\nabla_\xi R_\xi=0$. Then
$M$ holds $R_\xi\phi S=S\phi R_\xi$ if and only if $M$ is locally
congruent to one of the following:

{\setlength{\baselineskip}{0mm}
\begin{enumerate}
\item[$({\rm I})$] in case that $M_n(c)=P_n\mathbb{C}$ with $\eta(H\xi)\neq
0$,
\begin{enumerate}
\item[$(A_1)$] a geodesic hypersphere of radius $r$, where $0<r<
\pi/2$ and $r\neq \pi/4$,
\item[$(A_2)$] a tube of radius $r$
over a totally geodesic $P_k\mathbb{C}(1 \leq k \leq n-2)$, where
$0<r< \pi/2$ and $r\neq \pi/4$;
\end{enumerate}
\item[$({\rm II})$] in case that $M_n(c)=H_n\mathbb{C}$,
\begin{enumerate}
\item[$(A_0)$] a horosphere,
\item[$(A_1)$] a geodesic hypersphere or a tube over a complex hyperbolic hyperplane $H_{n-1} \mathbb{C}$,
\item[$(A_2)$] a tube over a totally geodesic $H_k\mathbb{C}(1 \leq k \leq
n-2)$.
\end{enumerate}
\end{enumerate}  }
\end{thm}
Proof. \underline{First step}. We prove $\beta =0$. \par

Suppose that $\beta\ne 0$. It follows from (4.22) that (3.10) is equivalent to
\begin{equation*}
(\rho_p \sigma_q -\sigma_p \rho_q)\phi_{pq}=0,
\end{equation*}
where $\rho_p =\alpha\lambda_p +c, \ \sigma_p =h\lambda_p -(\lambda_p)^2
+K$. Therefore if $\phi_{qp}\ne 0$, then we have
\begin{equation}
(\lambda_p-\lambda_q)\{-ch +\alpha\lambda_p \lambda_q +c(\lambda_p
 +\lambda_q) +\alpha K\} =0.
\end{equation}

Here we assert that if $\phi_{pq}\ne 0$, then $\lambda_p =\lambda_q$. To
 prove this, we assume that there exist indices $p$ and $q$ such that
\begin{equation*}
\phi_{pq}\ne 0, \quad \lambda_p -\lambda_q\ne 0.
\end{equation*}

First we prepare three Lemmas.
\begin{lem}
$(K\alpha^2 -c\alpha h)_1 =0$.
\end{lem}
Proof. From (5.1) we have
\begin{equation*}
(\alpha^2 K -\alpha hc) +(\alpha\lambda_p)(\alpha\lambda_q)
 +c(\alpha\lambda_p +\alpha\lambda_q) =0.
\end{equation*}
Lemma 3 follows from this and (4.26). \hfill$\Box$\vspace{0.2cm}\par

\begin{lem}
$4n\alpha\alpha_1 -(\alpha\gamma)_1 =0$.
\end{lem}
Proof. From (4.26) we have $(\alpha\sum_p \lambda_p)_1
 =0$. Combining this equation with $h=\alpha +\gamma +\delta
 +\sum_p\lambda_p$, we have 
\begin{equation*}
(\alpha(h -\alpha -\gamma -\delta))_1 =0.
\end{equation*}
Eliminate $h$ from this and Lemma 3. \hfill$\Box$\vspace{0.2cm}\par
\begin{lem}
$(\gamma -\delta -2n\alpha)\alpha_1 =0$ and $(\gamma -\delta
 -2n\alpha)\beta_1 =0$.
\end{lem}
Proof. From (4.24) we have $2\beta\beta_1 -(\alpha\gamma)_1 =0$. Hence
 it follows from Lemma 4 that
\begin{equation}
2n\alpha\alpha_1 -\beta\beta_1 =0.
\end{equation}
On the other hand, by (4.32) and (4.34) we have $(\gamma-\delta)\delta_1
 -\beta\delta_2 =0$, and therefore $(\gamma -\delta)\alpha_1
 -\beta\alpha_2 =0$. 
Thus Lemma 5 follows from (4.27) and (5.2). \hfill$\Box$\vspace{0.2cm}\par

We need to consider four cases.\vspace{0.2cm}\\
 Case I : Suppose that $\alpha_1\ne 0$ and $X_1 =0$. 
Owing to Lemma 5, we have $\gamma -\delta -2n\alpha =0$. Seeing $\theta_3$
 of \ $d$ \ of this equation and making use of (4.29), (4.30) and (4.33),
 we have 
\begin{equation}
2n\alpha^2(2n\alpha^2 -\delta^2 +2nc) +\beta^2\{3\delta^2 +(6n +4)c
 -2n\alpha^2\} =0.
\end{equation}
Seeing $\theta_1$ of \ $d$ \ of (5.3) and taking account of (3.2) and (5.2),
we have 
\begin{equation}
4n^2\alpha^4 +2n\alpha^2\{3\delta^2 +(8n +4)c\} -\beta^2(3\delta^2
 +2n\alpha^2) =0.
\end{equation}
Eliminating $\beta$ from (5.3) and (5.4), we have a polynomial of degree four
 with respect to $\delta$ containing the term $12n\alpha^2\delta^4\ne
 0$. This shows that $\delta$ is constant since $\alpha\delta +c =0$,
 which contradicts the assumption of Case I.\vspace{0.4cm}\\  
 Case II : Suppose that $\alpha_1\ne 0$ and $X_1 \ne 0$. By (4.25) we have
\begin{equation*}
\beta^2 -\alpha\gamma -c =0.
\end{equation*}Then from (4.39) we have
\begin{equation*}
(-\lambda_p\lambda_q +2c)(\lambda_p +\lambda_q) +2(\alpha
 +\gamma)\lambda_p\lambda_q -2c\gamma =0.
\end{equation*}
Multiply above equation by $\alpha^3$ and see $\theta_1$ of \ $d$ \ of this
equation. 
Then, from Lemma 4 and (4.26) we have
\begin{equation*}
c(\alpha\lambda_p +\alpha\lambda_q -\alpha\gamma)
 +(2n+1)(\alpha\lambda_p)(\alpha\lambda_q) -2cn\alpha^2 =0.
\end{equation*}
Again, seeing $\theta_1$ of \ $d$ \ of above equation, we have $cn\alpha\alpha_1
 =0$, which is a contradiction.\vspace{0.4cm}\\
 Case III : Suppose that $\alpha_1 = 0$ and $\beta^2 -\alpha\gamma -c
 \ne 0$. From (4.24), (4.25), (4.27), (4.28), (4.32) and (4.34) we have 
\begin{equation}
\delta_1 =\alpha_2 =\delta_2 =X_3 =\beta_1 =\gamma_1 =\beta_2 =X_1 =0.
\end{equation}
Seeing $\theta_2 \wedge \theta_3$ of \ $d$ \ of $\theta_{23}$ we have
$\beta_3 -2\beta^2 =\gamma\delta +2c$, which, together with (4.31) and
(5.5), imply
\begin{equation*}
\alpha\delta -\gamma\delta -\beta^2 =\gamma\delta +c.
\end{equation*}
Substituting of (4.14) and (5.5) into (4.29) we have
\begin{equation}
\alpha\gamma -\gamma\delta +\beta ^2 =-c.
\end{equation}
Eliminating $\beta$ from above two equations, we have
\begin{equation}
\alpha\delta -3\gamma\delta +\alpha\gamma =0.
\end{equation}
Seeing $\theta_2$ of \ $d$ \ of (5.6) and (5.7), we have
$(\alpha -\delta)\gamma_2 =0$ and $(\alpha -3\delta)\gamma_2
 =0$. Hence we have $\gamma_2 =0$.\par

Now put $F=\alpha,\beta,\gamma$ and $i=1,j=2$ in Lemma 1. Then,
 we have
\begin{equation*}
\alpha_3\gamma =\beta_3\gamma =\gamma_3\gamma =0.
\end{equation*}
If $\gamma\ne 0$, then from (4.14) and (4.33) we have a contradicton.
Thus $\gamma =0$, which contradicts (5.7). \vspace{0.4cm}\\ 
 Case IV : Suppose that 
\begin{align}
&\alpha_1 = 0,\\
&\beta^2 -\alpha\gamma -c = 0.
\end{align}
Seeing $\theta_2$ of \ $d$ \ of (5.9), we have 
\begin{equation}
(\beta^2 -\alpha\gamma)_3 =2\beta\beta_3 -\gamma\alpha_3 -\alpha\gamma_3
 =0.
\end{equation}
From (4.29)--(4.31), (4.33) and (5.9) we have the following: 
\begin{align}
&-\delta\gamma -\beta X_2 +(\gamma -\delta)X_1=0,\\
&\alpha_3 +3\beta\delta -\alpha\beta +\beta X_1 =0,\\
&\beta_3 +(\gamma -\delta)X_1 +\gamma\delta -\alpha\gamma -c =0,\\
&\gamma_3 -2\beta\delta +(\gamma -\delta)X_2 +\beta\gamma =0.
\end{align}\par
Substituting of (5.12)--(5.14) into (5.10) we have
\begin{equation*}
(\delta -\gamma)(X_1 -4\alpha)=0,
\end{equation*}
by virtue of (5.11). If $\delta =\gamma$, then by (5.9) we have a
 contradiction. Thus
\begin{equation}
X_1 =4\alpha.
\end{equation}
Substituting of this equation into (5.11)--(5.13) we have
\begin{align}
&\beta X_2 =4\alpha(\gamma-\delta)-\delta\gamma,\\
&\alpha_3 +3\beta\delta +3\alpha\beta =0,\\
&\beta_3 +3\alpha\gamma -3\alpha\delta +\gamma\delta =0.
\end{align}
It follows from (4.33), (5.9) and (5.16) that
\begin{equation}
\alpha\gamma_3 +\beta(3\alpha\gamma -6\alpha\delta -\gamma\delta) =0.
\end{equation}From (4.32), (5.2) and (5.8) we
 have $X_3 =0$ and $\beta_1 =0$ and therefore $\alpha_2 =\delta_2 =0$
 because of (4.27). 
Hence, seeing $\theta_1$ of \ $d$ \ of (5.9), we have $\gamma_1 =0$, and so
 $\beta_2 =0$. 


Now put $F=\alpha$ and $\beta$ in Lemma 1. Then we have 
\begin{equation*}
\alpha_3(\gamma +X_1) =0, \quad  \beta_3(\gamma +X_1) =0.
\end{equation*}
If $\gamma +X_1\ne 0$, then we have $\alpha_3 =\beta_3 =0$. It follow
 from (4.23) and (4.35) that $\alpha,\beta$ and $\delta$ are constant
 and that $\alpha_i =\beta_i =0$ for $i=1,2$. Furthermore, by (5.9) we see that
 $\gamma$ is constant. Thus from (5.17)-(5.19) we have  
\begin{align*}
&\alpha +\delta =0,\\
&3\alpha\gamma -3\alpha\delta +\gamma\delta =0,\\
&3\alpha\gamma -6\alpha\delta -\gamma\delta =0.
\end{align*}
Hence, by (3.2) and (5.9) we have $\alpha^2 -c=0$ and $2\beta^2 +c
 =0$, which is a contradiction. Therefore $X_1 =-\gamma$, which, together
 with (5.15), implies $\gamma =-X_1 =-4\alpha$. Thus it follows from (5.17) that
$\gamma_3 =-4\alpha_3=12\beta(\delta +\alpha)$. Hence from (5.19) we
 have a contradiction $\alpha\delta =0$.\vspace{0.2cm}\par 

Consequently, for all $p,q$ such that $\phi_{pq}\ne 0$, we have
 $\lambda_p =\lambda_q$. We take $p,q$ such that $\phi_{pq} \ne 0$ and put
 $\lambda :=\lambda_p =\lambda_q$. Then by (4.39) we have 
\begin{equation}
\beta^2 \lambda -(\lambda -\gamma)(\lambda^2 -\alpha\lambda -c) =0.
\end{equation}From $(1pq)$ and (4.38) we have
\begin{equation}
\lambda^2 +\alpha\lambda +c =0.
\end{equation}
Seeing $\theta_1$ of \ $d$ \ of above equation, we have $\lambda_1 =0$.
This and (4.26) imply $\alpha_1=\delta_1 =0$.  Thus it follows from
 (4.32), (4.34) and (4.27) that $X_3 =\alpha_2 =\delta_2 =\beta_1 =0$.
Substituting (5.21) into (5.20) we have
 $\beta^2 -(\lambda-\gamma)\lambda =0$. Hence seeing $\theta_1$ of \ $d$ \
 of this
 equation, we have $\gamma_1 =0$. Thus, from (4.28) we have $\beta_2
 =0$. On the other hand, from (5.21) we have 
\begin{equation*}
\beta^2 -\alpha\gamma -c =\beta^2 +\lambda^2\ne 0.
\end{equation*}
Thus by (4.25) we have $X_1 =0$. Therefore by the same argument as that
 in Case III, we have a contradiction. Consequantly we proved $\beta
 =0$. \vspace{0.2cm}\\ 

\noindent\underline{Second step}. From $(1pq)$ we have
\begin{equation*}
\alpha(2\lambda_p\lambda_q -\alpha\lambda_q -\alpha\lambda_p -2c)\phi_{qp} =0.
\end{equation*}
It follows from $(q1p)$ and (3.7) that
\begin{equation*}
(\lambda_p\lambda_q -\alpha\lambda_p -c)\phi_{pq} =0.
\end{equation*}From above two equations we have
\begin{equation}
\alpha^2 (\lambda_q -\lambda_p)\phi_{pq} =0.
\end{equation}
The equations $(312)$ and $(213)$ imply
\begin{equation}
\begin{split}
&(\alpha -\delta)\gamma +(\gamma -\delta)X_1 =0,\\
&(\gamma -\alpha)\delta +(\gamma -\delta)X_1 =0.
\end{split}
\end{equation}From above two equations we have
\begin{equation*}
(\gamma -\delta)(\alpha +2X_1) =0.
\end{equation*}
If $\gamma\ne \delta$, then $\alpha +2X_1 =0$ and so $X_1 \ne 0$.  Then
 by (3.4) we have a contradiction
\begin{equation*}
\alpha\gamma +c =\alpha(\gamma -\delta) =0.
\end{equation*}
Thus we have $\gamma =\delta$. From this and (5.23) we have $\alpha
 =\gamma$. Hence it follows from (5.22) that $H\phi =\phi H$.
Owing to Okumura's work or Montiel and Romero's work stated in the Introduction, we
 complete the proof of our Main Theorem. \hfill$\Box$\vspace{0.4cm}\par

\bigskip\par

\vspace{0.5cm} {\setlength{\baselineskip}{5mm}
\begin{quote}
U-Hang Ki \\
Department of Mathematics\\
Kyungpook National University \\
Daegu 702-701, Korea\\
E-mail address : uhangki2005@yahoo.co.kr \\

Hiroyuki Kurihara \\
Department of Computer and Media Science \\
Saitama Junior College \\
Hanasaki-ebashi, Kazo, Saitama 347-8503, Japan \\
E-mail address : kurihara@sjc.ac.jp \\

Ryoichi Takagi \\
Department of Mathematics and Informatics \\
Chiba University \\
Chiba 263-8522, Japan \\
E-mail address : takagi@math.s.chiba-u.ac.jp
\end{quote}}

\end{document}